\documentclass[reqno]{amsart}
\usepackage{hyperref}
\usepackage{amsmath,amssymb} 
\usepackage{graphicx}
\DeclareGraphicsExtensions{.pdf,.png,.jpg}

\usepackage[T2A]{fontenc}
\usepackage[utf8]{inputenc}
\usepackage{textcomp}

\usepackage{makecell}

\allowdisplaybreaks

\theoremstyle{plain}

\theoremstyle{definition}

\theoremstyle{remark}

\numberwithin{equation}{section}

\begin{document}
	\title[ Comparative analysis of Jacobi and Gauss-Seidel iterative methods
	]
	{Comparative analysis of Jacobi and Gauss-Seidel iterative methods
	}
	
	\author[P. Khrapov]{Pavel Khrapov}
	\address{Pavel Khrapov \\ Department of Mathematics
		\\ Bauman Moscow State Technical University \\  ul. Baumanskaya 2-ya, 5/1, Moscow \\ 105005, Moscow,  Russian Federation}  
	\email{khrapov@bmstu.ru , pvkhrapov@gmail.com }
	
	\author[N. Volkov]{Nikita Volkov}
	\address{Nikita Volkov \\ Department of Mathematics
		\\ Bauman Moscow State Technical University \\  ul. Baumanskaya 2-ya, 5/1, Moscow \\ 105005, Moscow,  Russian Federation}  
	\email{volkovns@student.bmstu.ru, nikita.volkov01@mail.ru}	
	
	
	
	\begin{abstract}
		The paper presents a comparative analysis of iterative numerical methods of Jacobi and Gauss-Seidel for solving systems of linear algebraic equations (SLAEs) with complex and real matrices. The ranges of convergence for both methods for SLAEs in two and three unknowns, as well as the interrelationships of these ranges are obtained. An algorithm for determining the convergence of methods for SLAEs using the complex analog of the Hurwitz criterion is constructed, the realization of this algorithm in Python in the case of SLAEs in three unknowns is given. A statistical comparison of the convergence of both methods for SLAEs with a real matrices and the number of unknowns from two to five is carried out.
	\end{abstract}
	
	\maketitle
	
	\section{Introduction} 
	
	\indent In the modern world, a large number of both applied and theoretical problems in various fields of science and technologies are reduced to the problem of finding exact solutions of various SLAEs or solutions that maximally approximate the exact ones, numerical methods for solving which have been developing over the years due to the huge number of areas of their application \cite{Bylina}, \cite{Nutzi}.\\
	\indent A special place in the theory of SLAEs' solutions is occupied by the simple iterative method, which is an alternative to direct methods of finding SLAEs' solutions. At the same time, based on the simple iterative method, new methods for solving SLAEs are being developed, which are an improved version of the classical method \cite{Saad}, \cite{Tarigan}, \cite{Gunawardena}.\\
	\indent Some of these, based on the simple iterative method, are the Jacobi and Gauss-Seidel iterative methods for solving SLAEs, the meaning of which is to allocate elements on, above and below the diagonal of the original SLAE's matrix as separate matrices and conduct the simple iterative method using them instead of the original, which often greatly simplifies
	the calculations \cite{Bagnara}, \cite{Salkuyeh}. Iterative Jacobi and Gauss-Seidel methods, also being classical iterative methods of solving SLAEs, have recently undergone various improvements, some of which are described for example in \cite{Salkuyeh}, \cite{Milaszewicz}, \cite{Sun}, \cite{Ahmadi}, \cite{Kohno}, \cite{Li1}, \cite{Shang}, \cite{Courtecuisse}, \cite{Tavakoli}. Nevertheless, many modern alternatives to the classical Jacobi and Gauss-Seidel methods are based on sufficient condition of their convergence to an exact solution in the case of diagonal predominance in the original SLAEs' matrices, without considering the cases without diagonal predominance when these methods can also converge to an exact solution, and are also described only for special types of matrices \cite{Koester}, \cite{Amodio}.\\
	\indent The convergence of iterations to an exact solution is one of the main problems, since, as a consequence of the classical simple iterative method, the Jacobi and Gauss-Seidel methods not always converge to an exact solution, and have convergence criteria following from a similar criterion for the simple iterative method \cite{Bagnara}. The search for ranges of convergence and the theoretical comparison of the effectiveness of the methods based on it is the main task of this work.\\
	\indent The convergence criteria obtained in \cite{Bagnara}, according to which the eigenvalues of
	the matrices in the method should be less than one in absolute value, are reduced to the problem of finding the roots of the algebraic polynomials of degree $n$ with complex coefficients inside the unit circle, various solutions of which are described for example in \cite{Chen}, \cite{Bharanedhar}, \cite{Zadorojniy}, \cite{Konvalina}, \cite{Korsakov}, \cite{Biberdorf}, \cite{Joyal}, \cite{Dehmer}, \cite{Frank}, and for polynomials of a special kind in \cite{Erdelyi}, \cite{Mercer}.\\ 
	\indent It can be solved by making a fractional linear transformation that translates the interior of the unit circle of the complex plane to the left half-plane and reduces it to the study of \textbf{stability} of the polynomial \cite{Zadorojniy}. In \cite{Zadorojniy} this problem is considered for polynomials with real coefficients of the second and third degree.\\
	\indent In this paper, a comparative analysis of two methods using the examples of SLAEs in two and three unknowns is carried out by considering the ranges of their convergence, which are obtained under the assumption that the boundary of each range is formed when at least one root of the corresponding equation has a unit absolute value, and the rest does not exceed one, and all points are contained inside the range, for which all roots have absolute values less than one.\\
	\indent There is described  the general convergence criteria for each method in paragraph 2.\\
	\indent In paragraph 3, the convergence ranges of the methods for SLAEs with complex coefficients in two unknowns are obtained, and the conclusion of their comparison is given: for the Jacobi method the convergence range and its boundary are found when substituting roots with absolute values less than one into the corresponding equation, and for the Gauss-Seidel method by directly solving the equation.\\
	\indent In paragraph 4, by a similar substitution of roots with absolute values not exceeding one, the boundary conditions of the methods for SLAEs with complex coefficients in three unknowns are obtained, and on their basis the convergence ranges in the real case are obtained, for which a comparative analysis is given.\\
	\indent  There is described a general convergence check method for SLAEs with complex matrix based on \cite{Zadorojniy} and \cite{Postnikov}, and a general comparison of both methods is made in paragraph 5 .\\
	\indent In paragraph 6, a statistical comparison of convergence of both methods for SLAE with real matrix is carried out using mathematical modeling.
	
	\section{Convergence conditions of the Jacobi and Gauss-Seidel methods} 
	\label{Conditions}
	
	\indent When solving a system of linear algebraic equations
	\begin{equation}\label{1}
		A\vec{x}=\vec{b}
	\end{equation}
	in accordance with the Jacobi method, the matrix $A$ of the original SLAE is represented as a sum:
	\[
	A = L + D + R
	\]
	\[
	\det A \neq 0
	\]
	where $L, D, R,$ are, respectively, the matrices with subdiagonal, diagonal, and overdiagonal  elements of matrix $A$, and then there is a system is obtained from the original SLAE~(\ref{1}):
	\[
	\vec{x} = -D^{-1}(L + R)\vec{x} + D^{-1}\vec{b}
	\]
	for which the simple iterative method converges if all roots of the equation
	\begin{equation}\label{2}
		\begin{vmatrix}
			\lambda a_{11}&a_{12}&...&a_{1n}\\
			a_{21}&\lambda a_{22}&...&a_{2n}\\
			...& ...& ...& ...\\
			a_{n1}& a_{n2}& ...& \lambda a_{nn}
		\end{vmatrix} = 0
	\end{equation}
	have absolute values less than one \cite{Bagnara}; $a_{ij}$ - elements of the original matrix $A$, $a_{ij} \in \mathbb{C}$.\\
	\indent Similarly, the Gauss-Seidel method transforms the original SLAE~(\ref{1}) to a system:
	\[
	\vec{x} = -(L + D)^{-1}R\vec{x} + (L + D)^{-1}\vec{b}
	\]
	for which the simple iterative method converges to an exact solution if all roots of equation
	\begin{equation}\label{3}
		\begin{vmatrix}
			\lambda a_{11}& a_{12}& ...& a_{1n}\\
			\lambda a_{21}& \lambda a_{22}& ...& a_{2n}\\
			...& ...& ...& ...\\
			\lambda a_{n1}& \lambda a_{n2}& ...& \lambda a_{nn}
		\end{vmatrix} = 0
	\end{equation}
	have absolute values less than one \cite{Bagnara}, $a_{ij} \in \mathbb{C}$.\\
	\indent When the dimension of the original SLAE is small, we can find the convergence ranges of the methods by directly solving the equations~(\ref{2}) and~(\ref{3}). Let us show this for the cases of SLAEs in two and three unknowns, which often arise in applied research.
	
	\section{System of linear algebraic equations in two unknowns}
	\label{Two unknowns}
	\subsection{Jacobi method}
	\indent The equation~(\ref{2}) has the form:
	\begin{equation}\label{4}
		\begin{vmatrix}
			\lambda a_{11}& a_{12}\\
			a_{21}& \lambda a_{22}
		\end{vmatrix}
		= \lambda ^2 a_{11} a_{22} - a_{12} a_{21} = 0,
	\end{equation}
	and for the convergence of the method it is necessary that its roots lie inside the unit circle. \\
	\indent In the general case $a_{ij}, \lambda_{1,2}\in\mathbb{C}$, and the system~(\ref{2}) is equivalent to equations
	\[
	(\lambda-r_{1}e^{i\varphi_{1}})(\lambda-r_{2}e^{i\varphi_{2}}) = 0
	\]
	\begin{equation}\label{5}
		\lambda^{2}+\lambda(-r_{1}e^{i\varphi_{1}}-r_{2}e^{i\varphi_{2}}) +r_{1}r_{2}e^{i\varphi_{1}}e^{i\varphi_{2}} = 0,
	\end{equation}
	 where $r_1e^{i\varphi_{1}}, r_2e^{i\varphi_{2}}$ are the roots of equation~(\ref{4}), $r_{1},r_{2}<1$.\\
	 \indent Comparing~(\ref{4}) and~(\ref{5}), we obtain the system ($a_{11}a_{22} \neq 0$, since $a_{11}, a_{22}$ are elements of the diagonal matrix $D$):
	 \begin{equation}\label{6}
	 	\begin{cases}
	 		r_{1}e^{i\varphi_{1}}+r_{2}e^{i\varphi_{2}}=0\\
	 		r_{1}r_{2}e^{i\varphi_{1}}e^{i\varphi_{2}}=-\frac{a_{12}a_{21}}{a_{11}a_{22}}
	 	\end{cases}
	 \end{equation}
 	from which follow:
 	\[
 	r_1=r_2
 	\]
 	\[
 	r_{1}r_{2}=r_{1}^2 = |\frac{a_{12}a_{21}}{a_{11}a_{22}}|<1
 	\]
 	\begin{equation}\label{7}
 		|a_{12}a_{21}|<|a_{11}a_{22}|
 	\end{equation}
 	\indent The condition~(\ref{7}) defines the convergence range for the Jacobi method in the general case: the absolute value of the product of the off-diagonal elements of the matrix $A$ of the system~(\ref{1}) must be less than the absolute value of the product of its diagonal elements for the method to converge in the case of an SLAE in two unknowns.\\
 	\subsection{Gauss-Seidel method}
 	\indent The equation~(\ref{3}) has the form:
 	\begin{equation}\label{8}
 		\begin{vmatrix}
 			\lambda a_{11}& a_{12}\\
 			\lambda a_{21}& \lambda a_{22}
 		\end{vmatrix}
 		= \lambda ^2 a_{11} a_{22} - \lambda a_{12} a_{21} = 0
 	\end{equation}
 	and its roots
 	\[
 	\lambda_1 = 0
 	\]
 	\[
 	\lambda_2 = \frac{a_{12} a_{21}}{a_{11} a_{22}}.
 	\]
 	must have an absolute value less than one ($a_{11}a_{22} \neq 0$, since $a_{11}, a_{22}$ are diagonal elements of the triangular matrix $L+D$).\\
 	\indent Since one of them is zero, only the second root is checked for the convergence condition, for which, in order for its absolute value to be less than one, it is necessary to fulfill the condition
 	\begin{equation}\label{9}
 		|a_{12}a_{21}| < |a_{11}a_{22}|
 	\end{equation}
 	\indent \textbf{Thus, both the Jacobi method and the Gauss-Seidel method for SLAEs in two unknowns have the same range of convergence}~(\ref{9})\textbf{.}\\

	\section{System of linear algebraic equations with three unknowns}
	\label{Three unknowns}
	\subsection{Jacobi method}
\begin{figure}[h]
			\begin{minipage}[h]{0.4\linewidth}
				\includegraphics[height=6.5cm, width=6.9cm]{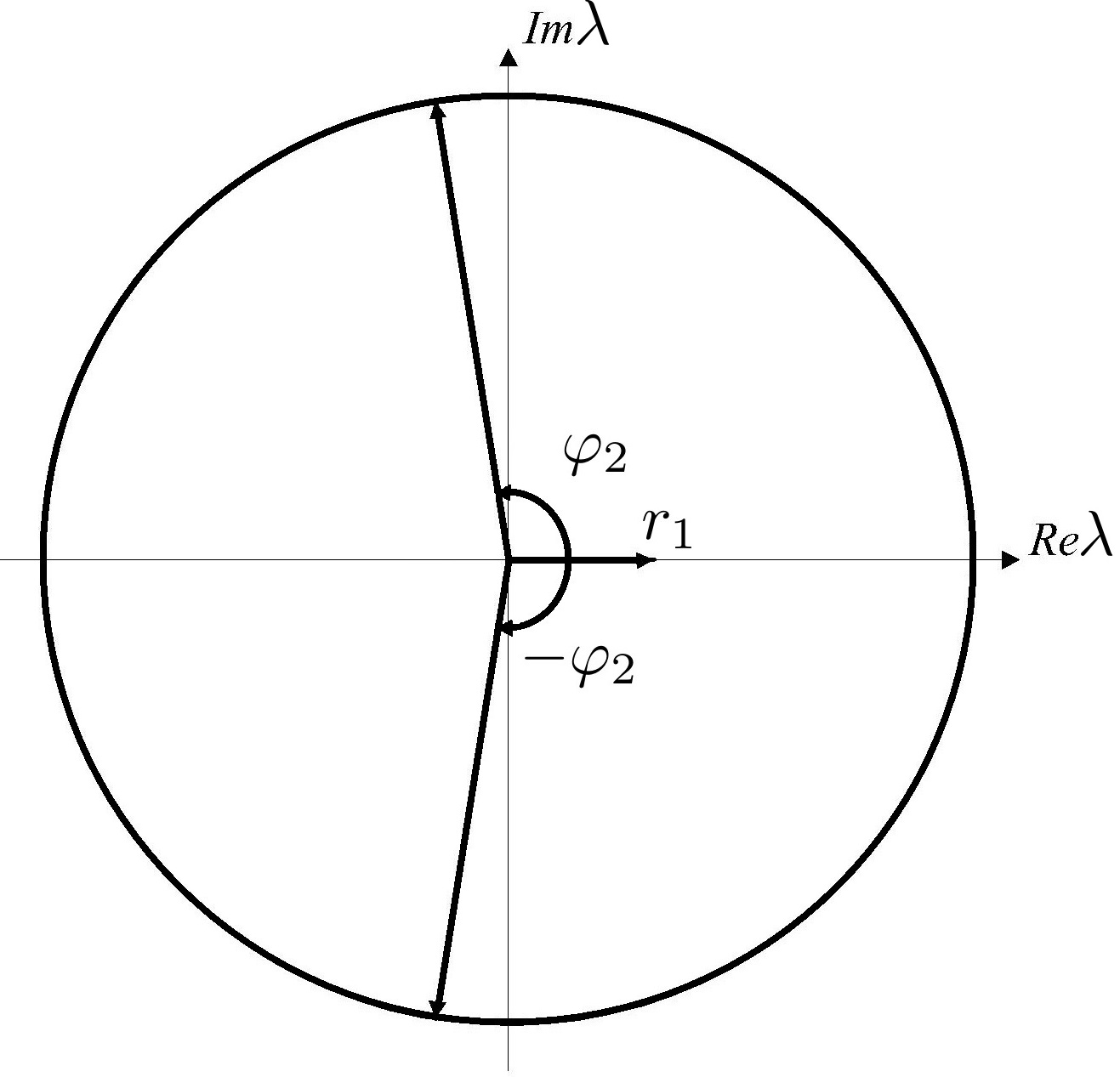}
		\caption{Location of roots of the equation~(\ref{11}) on the complex plane before rotation}
				\label{p4}
			\end{minipage}
			\hfill
			\begin{minipage}[h]{0.4\linewidth}
				\includegraphics[height=6.5cm, width=6.9cm]{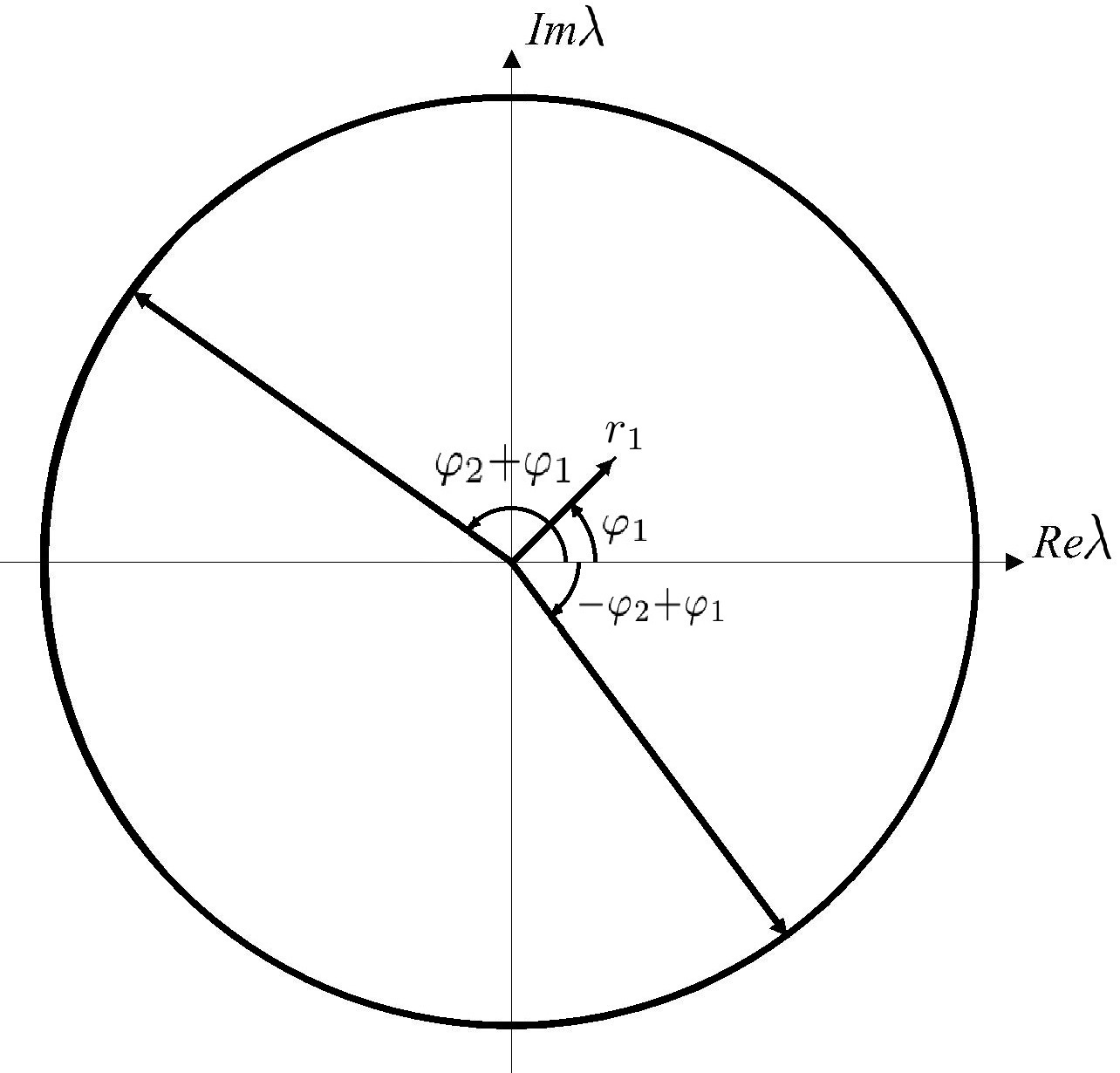}
		\caption{Location of roots of the equation~(\ref{11}) on the complex plane after rotation}
				\label{p5}
			\end{minipage}
	\end{figure}

	\indent The equation~(\ref{2}) has the form:
	\begin{equation}\label{10}
		\lambda ^3 a_{11}a_{22}a_{33} - \lambda (a_{13}a_{22}a_{31} + a_{23}a_{11}a_{32} + a_{12}a_{33}a_{21}) + (a_{13}a_{32}a_{21} + a_{12}a_{23}a_{31}) = 0
	\end{equation}
	 and for convergence of the method it is necessary that all its roots lie inside the unit circle.\\
	 \indent Let's divide it by $a_{11}a_{22}a_{33}$ (there are no zero elements on the diagonal of matrix $A$, since matrix $D$ must have an inverse):
	 \[
	 \lambda ^3 + \frac{-a_{13} a_{22} a_{31} - a_{23} a_{11} a_{32} - a_{12} a_{33} a_{21}}{a_{11} a_{22} a_{33}}\lambda + \frac{a_{13} a_{32} a_{21} + a_{12} a_{23} a_{31}}{a_{11}a_{22}a_{33}} = 0.
	 \]
	 \indent Denoting
	 \[
	 p = \frac{-a_{13} a_{22} a_{31} - a_{23} a_{11} a_{32} - a_{12} a_{33} a_{21}}{a_{11} a_{22} a_{33}}
	 \]
	 \[
	 q = \frac{a_{13} a_{32} a_{21} + a_{12} a_{23} a_{31}}{a_{11}a_{22}a_{33}},
	 \]
	 we obtain the canonical cubic equation:
	 \begin{equation}\label{11}
	 	\lambda ^3 + p\lambda + q = 0,
	 \end{equation}
 	\indent In the general case its coefficients and roots are complex: $p, q, \lambda_{1,2,3}\in\mathbb{C}$.\\
 	\indent We find the \textbf{range of convergence} of the method expressed in terms of $p, q \in\mathbb{C}$ by obtaining the equations of its boundaries and combining them.\\
 	\indent We obtain the equations of the boundaries under the assumption that there is at least one root of the equation~(\ref{11}) on boundaries, the absolute value of which is equal to one, and the interior points of the range are those in which the absolute value of each root is less than one. The boundary is not included in the convergence range, since at least one of the roots has a unit absolute value on it, which contradicts the convergence condition \cite{Bagnara}. Consider several cases.\\
 	\indent Find the equation for the first boundary of the convergence range: let one of the roots of the equation~(\ref{11}) have a unit absolute value, and the other two roots have an absolute value not exceeding one:
 	\begin{equation}\label{12}
 		(\lambda-e^{i\varphi_{1}})(\lambda-r_{2}e^{i\varphi_{2}})(\lambda-r_{3}e^{i\varphi_{3}})=0
 	\end{equation}
 	\[
 	\lambda^3+\lambda^2(-e^{i\varphi_{1}}-r_{2}e^{i\varphi_{2}}-r_{3}e^{i\varphi_{3}}) +\lambda(r_{2}e^{i\varphi_{1}}e^{i\varphi_{2}} + r_{3}e^{i\varphi_{1}}e^{i\varphi_{3}} + r_{2}r_{3}e^{i\varphi_{2}}e^{i\varphi_{3}}) -   r_{2} r_{3}e^{i\varphi_{1}}e^{i\varphi_{2}}e^{i\varphi_{3}} = 0
 	\]
 	\[
 	r_{2},r_{3}\leq1
 	\]
 	\indent Comparing~(\ref{11}) and~(\ref{12}), we obtain a system of equations for the first boundary of the convergence range:
 	\begin{equation}\label{13}
 		\begin{cases}
 			e^{i\varphi_{1}} + r_{2}e^{i\varphi_{2}} + r_{3}e^{i\varphi_{3}} = 0\\
 			r_{2}e^{i\varphi_{1}}e^{i\varphi_{2}} + r_{3}e^{i\varphi_{1}}e^{i\varphi_{3}} + r_{2}r_{3}e^{i\varphi_{2}}e^{i\varphi_{3}} = p\\
 			-   r_{2} r_{3}e^{i\varphi_{1}}e^{i\varphi_{2}}e^{i\varphi_{3}} = q
 		\end{cases}
 	\end{equation}
 	from which follow:
 	\begin{equation}\label{14}
 		|q|\leq 1
 	\end{equation}
 	\begin{equation}\label{15}
 		\arg(q) = \pi + \varphi_{1} + \varphi_{2} + \varphi_{3}
 	\end{equation}
 	\begin{equation}\label{16}
 		r_{3}e^{i\varphi_{3}} = - r_{2}e^{i\varphi_{2}} - e^{i\varphi_{1}} 
 	\end{equation}
	\indent The first equation of the system~(\ref{13}) has a geometric interpretation (fig.~\ref{p1}):
	\begin{figure}[h]
		\center
		\includegraphics[height=6cm, width=6cm]{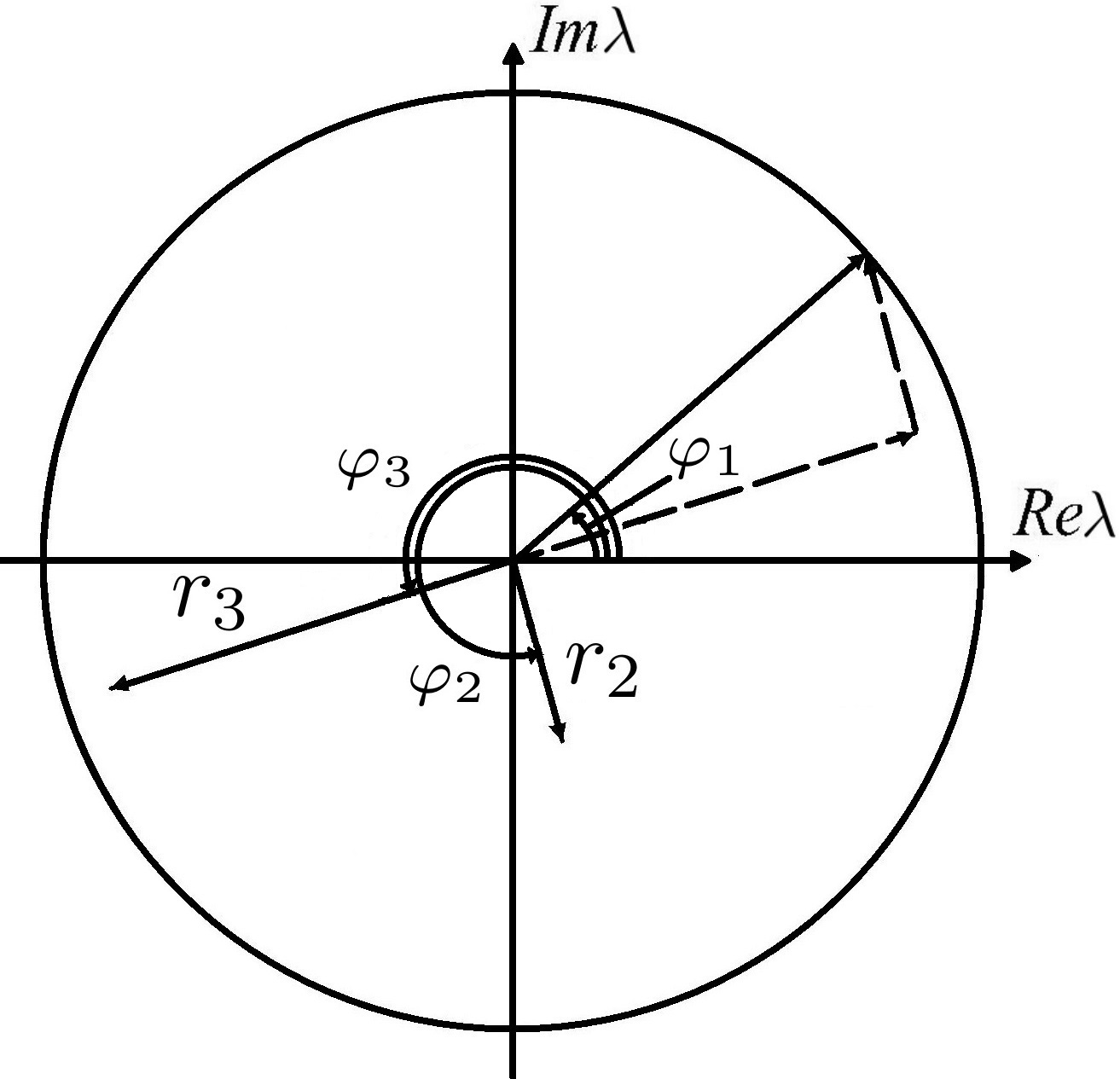}
		\caption{Geometric interpretation of the first equation of the system~(\ref{13})}
		\label{p1}
	\end{figure}\\
	\indent  Substituting~(\ref{16}) into the second and third equations of the system~(\ref{13}), we have expressions:
	\[
	p = -e^{i\varphi_{1}}e^{i\varphi_{1}} - r_{2}e^{i\varphi_{1}}e^{i\varphi_{2}} -r_{2}^{2}e^{i\varphi_{2}}e^{i\varphi_{2}}
	\]
	\[
	q = r_{2}e^{i\varphi_{1}}e^{i\varphi_{1}}e^{i\varphi_{2}} + r_{2}^{2}e^{i\varphi_{1}}e^{i\varphi_{2}}e^{i\varphi_{2}}
	\]
	comparing which, we obtain the equation of the first boundary of the convergence range of the Jacobi method in the general case:
	\begin{equation}\label{17}
		p = -qe^{-i\varphi_{1}} - e^{2i\varphi_{1}}
	\end{equation}
	 from which follow the relationships of absolute value and argument for $p$ and $q$:
	\begin{equation}\label{18}
		|p|=\sqrt{r_{q}^2 + 1 + 2r_{q}\cos(\varphi_{q} - 3\varphi_{1})}
	\end{equation}
	\begin{equation}\label{19}
		\arg(p) = \arctan(\frac{r_{q}\sin(\varphi_{q} - \varphi_{1}) +\sin(2\varphi_{1})}{r_{q}\cos(\varphi_{q} - \varphi_{1}) + \cos(2\varphi_{1})})
	\end{equation}
	\indent The relationships~(\ref{18}) and~(\ref{19}) show that the absolute value and argument for $p$ depend on three parameters - the absolute value $r_q$, the argument $\varphi_{1}$, and the argument $\varphi_{q}$, the last of which depends not only on $\varphi_{1}$~(\ref{15}), so we will take the argument $\varphi_{1}$ as a parameter to visualize the absolute value and argument for $p$.\\
	\indent Let us take for example $\varphi_{1} = 0$ and $\varphi_{1} = \pi$ (we use these parameters for further visualization of the special case of SLAE with real matrix when at least one of the roots of the equation~(\ref{11}) is real): at $\varphi_{1} = 0$, according to geometrical considerations (fig.~\ref{p1}) and~(\ref{15})
	\begin{equation}\label{20}
		\varphi_{q} = \varphi_{2} + \varphi_{3} + \pi \in [-\pi, -\frac{\pi}{2}] \cup [\frac{\pi}{2}, \pi]
	\end{equation}
	\begin{equation}\label{21}
		p = -q-1
	\end{equation}
	\indent Similarly, when $\varphi_{1} = \pi$
	\begin{equation}\label{22}
		\varphi_{q} = \varphi_{2} + \varphi_{3} + \pi + \pi  \in [-\frac{\pi}{2}, \frac{\pi}{2}]
	\end{equation}
	\begin{equation}\label{23}
		p = q-1
	\end{equation}
	\indent The picture of the absolute value~(\ref{18}) in this case is as follows (fig.~\ref{p2}):
	\begin{figure}[h]
		\center
		\includegraphics[height=7cm, width=9.9cm]{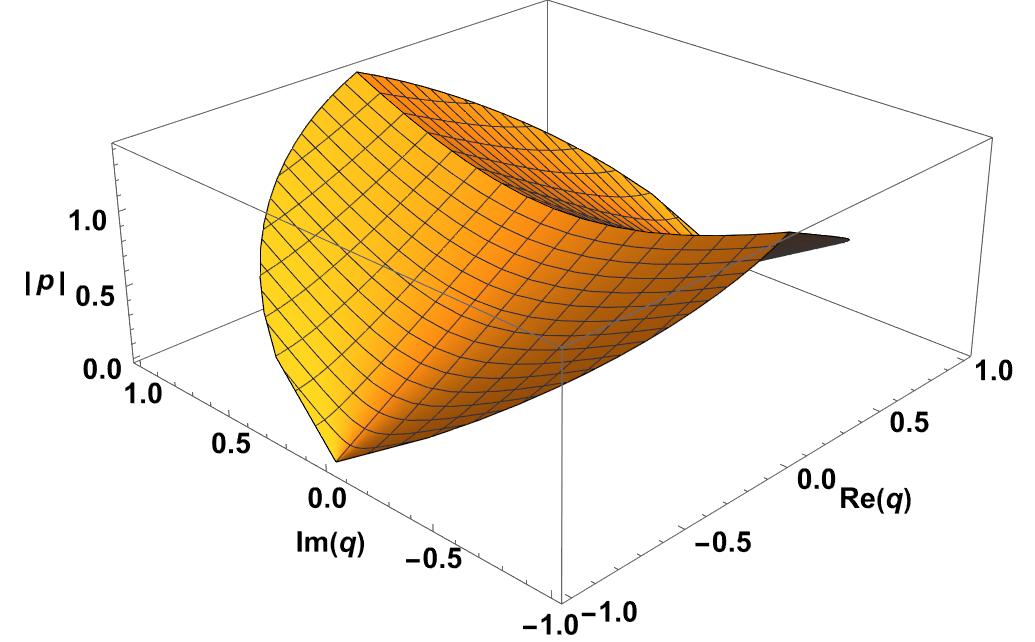}
		\caption{Dependence of the absolute value $|p|$ at the boundary~(\ref{17}) at $\varphi_{1} = 0$ ($\operatorname{Re} (q) \leq 0$) and $\varphi_{1} = \pi$ ($\operatorname{Re} (q) \geq 0$)}
		\label{p2}
	\end{figure}\\
	\indent To find the second boundary of the convergence range, consider the case when two roots of the equation~(\ref{11}) on the complex plane have axial symmetry with respect to the line passing through the vector of the third root (the case when two roots are complex-conjugate and the third is real is a special case of this case), and the roots located symmetrically have a unit absolute value, and the third root has an absolute value not exceeding one. \\
	\indent This case can be considered as a rotation of the system of vectors of roots of the equation on the complex plane from the zero angle by the angle $\varphi_{1}$, which is the argument of the first root: taking into account that before the rotation by the angle $\varphi_{1}$ one root was real, and the other two roots were complex-conjugate with arguments $\varphi_{2}$ and $-\varphi_{2}$ respectively, after the rotation the picture on the complex plane will be as follows (fig.~\ref{p4}, fig.~\ref{p5}):\\
	\begin{equation}\label{24}
		(\lambda-r_{1}e^{i\varphi_{1}})(\lambda-e^{i(\varphi_{2} + \varphi_{1})})(\lambda-e^{i(-\varphi_{2} + \varphi_{1})})=0
	\end{equation}
	\[r_1 \leq 1\]
	\indent Opening the brackets and comparing~(\ref{24}) with equation~(\ref{11}), we obtain the system for the second boundary of the convergence range:
	\begin{equation}\label{25}
		\begin{cases}
			e^{i(\varphi_{2} + \varphi_{1})} + e^{i(-\varphi_{2} + \varphi_{1})}  + r_{1}e^{i\varphi_{1}}= 0\\
			e^{2i\varphi_{1}} + r_{1}e^{(\varphi_{1} + \varphi_{2})}e^{i\varphi_{1}} + r_{1}e^{(\varphi_{1} - \varphi_{2})}e^{i\varphi_{1}} = p\\
			-r_{1}e^{3i\varphi_{1}} = q
		\end{cases}
	\end{equation}
	from which follow:
	\begin{equation}\label{26}
		|q|\leq 1
	\end{equation}
	\begin{equation}\label{27}
		r_{1}e^{i\varphi_{1}} = -e^{i(\varphi_{2} + \varphi_{1})} - e^{i(-\varphi_{2} + \varphi_{1})}
	\end{equation}
	\indent Let's substitute~(\ref{27}) into the second and third equations of the system~(\ref{25}):
	\begin{equation}\label{28}
		p = -e^{2i(\varphi_{1} + \varphi_{2})} - e^{2i\varphi_{1}} - e^{2i(\varphi_{1} - \varphi_{2})}
	\end{equation}
	\[
	q = e^{2i\varphi_{1}}(e^{i(\varphi_{1} + \varphi_{2})} + e^{i(\varphi_{1} - \varphi_{2})})
	\]
	\begin{equation}\label{29}
		qe^{-2i\varphi_{1}} = e^{i(\varphi_{1} + \varphi_{2})} + e^{i(\varphi_{1} - \varphi_{2})}
	\end{equation}
	\indent Comparing the square of the expression~(\ref{29}) and the expression~(\ref{28}), we obtain the equation of the second boundary of the convergence range of the Jacobi method in general case:
	\begin{equation}\label{30}
		p = -q^2e^{-4i\varphi_{1}} + e^{2i\varphi_{1}}
	\end{equation}
	\begin{equation}\label{31}
		\arg(q) = 3\varphi_{1} + \pi
	\end{equation}
	from which we find the relationships of absolute value and argument for $p$ and $q$ subject to the condition~(\ref{31}):
	\begin{equation}\label{32}
		|p| = 1-r_{1}^2 = 1 - |q|^2
	\end{equation}
	\begin{equation}\label{33}
		\arg(p) = 2\varphi_{1} = \frac{2}{3} \arg(q) - \frac{2}{3} \pi
	\end{equation}
	\indent When $\varphi_{1} \in [-\pi, \pi]$, the picture for the absolute value~(\ref{32}) is as follows (fig.~\ref{p8}):
	\begin{figure}[h]
		\center
		\includegraphics[height=7cm, width=9.9cm]{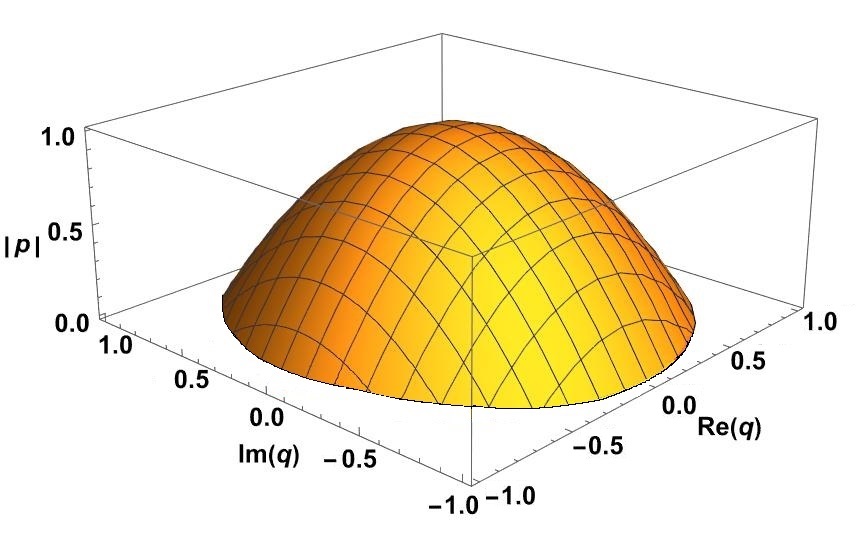}
		\caption{Dependence of the absolute value $|p|$ on $q$}
		\label{p8}
	\end{figure}\\
	\indent  It follows from~(\ref{31}),~(\ref{33}) that when $\varphi_{1} = \pm \pi$, arguments  
	\[
	\arg(q)=\arg(p)=0
	\] 
	are in the same phase corresponding to the real case, and when $\varphi_{1} = 0$ 
	\[
	\arg(q)= \pi
	\]
	\[
	\arg(p)=0
	\]
	the arguments are in antiphase corresponding also to the real case ($p$ is positive, $q$ is negative).\\
	\indent In addition, for the first bound~(\ref{17}), it follows from~(\ref{19}) that for $\varphi_{q} \in \{0, \pm \pi\}$ and $\varphi_{1} \in \{0, \pm \pi\}$ , as in the case of the second bound:
	\[
	\arg(p) \in \{0, \pm \pi\}
	\]
	\indent Hence, when at least one of the roots of the equation~(\ref{11}) is real, the cross section of the convergence range boundaries shown in the figures (fig.~\ref{p2}, fig. ~\ref{p8}), at $\operatorname{Im}(q) = 0$, can be combined into one general boundary of the convergence range of the Jacobi method in the real case (fig.~\ref{p9}) (by getting rid of absolute values, part of the boundaries, according to~(\ref{21}) and~(\ref{23}), moves to the area of negative values of $p$):
	\begin{figure}[h]
		\center
		\includegraphics[height=6cm, width=7.4cm]{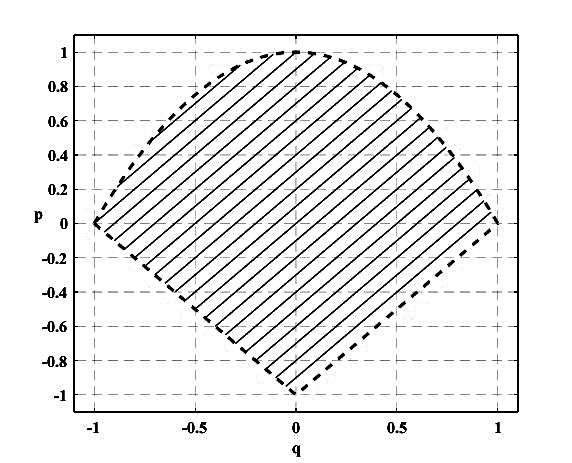}
		\caption{Boundary of the convergence range of the Jacobi method for the real case}
		\label{p9}
	\end{figure}\\
	\indent In particular, when $\varphi_{1} \in \{0, \pm \pi\}$, the equation~(\ref{30}) takes the form of a parabola:
	\[
	p = q^2 + 1
	\]
	\indent Thus, the boundary of the convergence range of the Jacobi method in the general complex case is the union of the sets of points satisfying the equations~(\ref{17}) and~(\ref{30}) and the condition:
	\[
	|q|\leq1
	\]
	\indent In the real case, when $a_{ij} \in A$; $p, q$ are real numbers, and the roots of equation~(\ref{11}) are either all real or two of them are complex-conjugate, equations~(\ref{17}) and~(\ref{30}) form the following boundary of the convergence range (fig.~\ref{p9}):
	\begin{equation}\label{34}
		\left[
		\begin{gathered}
			p= -q-1\\
			p= q - 1\\
			p = -q^{2} + 1\\
			-1 \leq q \leq 1
		\end{gathered}
		\right.
	\end{equation}
	\indent The boundary, due to the above assumptions that it contains at least one root of the equation~(\ref{11}) whose absolute value is equal to one, does not belong to the convergence range. The convergence range consists of the set of points bounded by the boundary~(\ref{34}), which does not belong to this region, so for the convergence range, given the conditions~(\ref{14}) and~(\ref{26}), it follows that:
	\[
	-1 < q < 1
	\]
	\indent Let's show that the area in the figure (fig.~\ref{p9}) can be filled completely:
	\[
	p = -\frac{a_{13}a_{31}}{a_{11}a_{33}} - \frac{a_{23}a_{32}}{a_{22}a_{33}} - \frac{a_{12}a_{21}}{a_{11}a_{22}} = 
	- \frac{a_{13}a_{32}a_{21}a_{22}a_{31}}{a_{11}a_{22}a_{33}a_{21}a_{32}} - \frac{a_{12}a_{23}a_{31}a_{11}a_{32}}{a_{11}a_{22}a_{33}a_{12}a_{31}} - \frac{a_{12}a_{21}}{a_{11}a_{22}}
	\]
	\indent Let's denote:
	\[
	x = \frac{a_{22}}{a_{21}}
	\]
	\[
	y = \frac{a_{11}}{a_{12}}
	\]
	\[
	t = \frac{a_{31}}{a_{32}}
	\]
	\[
	p = - \frac{a_{13}a_{32}a_{21}}{a_{11}a_{22}a_{33}}xt - \frac{a_{12}a_{23}a_{31}}{a_{11}a_{22}a_{33}}\frac{y}{t} - \frac{1}{xy}
	\]
	\indent Assuming that:
	\begin{equation}\label{35}
		xt = \frac{y}{t} = a
	\end{equation}
	\[
	t\neq0
	\]
	\[
	xy = a^2
	\]
	\begin{equation}\label{36}
		p = -aq - \frac{1}{a^2}
	\end{equation}
	\indent $x, y, t$ are independent of each other, so it is always possible to choose the coefficients of the matrix $A$ of the system~(\ref{1}) such that the condition~(\ref{35}) is satisfied and the line~(\ref{36}) is obtained. At the same time, $p$ and $q$ depend on three more parameters on which $x, y, t$ do not depend, so it is possible to choose $p$ and $q$ such that they lie in the convergence range. Thus, we can construct any number of lines of the form~(\ref{36}), some of whose points lie inside the convergence range labeled in the figure (Fig.~\ref{p9}). The set of such lines completely intersects the convergence range. Accordingly, it is always possible to find a SLAE for which $p$ and $q$ lie within the convergence range of the Jacobi method.\\
	\subsection{Gauss-Seidel method}
	\indent The equation~(\ref{3}) has the form:
	\begin{equation}\label{37}
		\lambda ^3 a_{11}a_{22}a_{33} + \lambda^2 (a_{21}a_{13}a_{32} - a_{13}a_{22}a_{31} - a_{32}a_{11}a_{23} -a_{21}a_{33}a_{12}) + \lambda a_{12}a_{23}a_{31} = 0
	\end{equation}
	and the method converges if all its roots lie inside the unit circle.\\
	\indent One of the roots of the equation~(\ref{37}) is zero, and the other two roots are found from the quadratic equation:
	\begin{equation}\label{38}
		\lambda ^2 a + \lambda d + b = 0
	\end{equation}
	where
	\[ 
	d = a_{21}a_{13}a_{32} - a_{13}a_{22}a_{31} - a_{32}a_{11}a_{23} -a_{21}a_{33}a_{12}
	\]
	\[
	a = a_{11}a_{22}a_{33}
	\]
	\[
	b=a_{12}a_{23}a_{31}
	\]
	\indent For the Gauss-Seidel method, we find the convergence range, given that within it all roots of the equation~(\ref{38}) have an absolute value not exceeding one.\\
	\indent In the general case $a, b, d, \lambda_{1,2} \in \mathbb{C}$ and the equation~(\ref{38}) is equivalent to equation:
	\[
	a\lambda^2-\lambda a(r_1 e^{i\varphi_{1}} + r_{2}e^{i\varphi_{2}}) + ar_1 r_{2}e^{i\varphi_{1}}e^{i\varphi_{2}} = 0	
	\]
	\[
	0< r_1 , r_{2} < 1
	\]
	comparing it with~(\ref{38}), obtain the system:
	\begin{equation}\label{39}
		\begin{cases}
			-ar_1 e^{i\varphi_{1}} - ar_{2}e^{i\varphi_{2}} = d \\
			ar_1r_{2}e^{i\varphi_{1}}e^{i\varphi_{2}} = b
		\end{cases}
	\end{equation}
	which defines the convergence range of the Gauss-Seidel method in the general case, and from which it follows that in the convergence range ($a \neq 0$, since $a_{ii}$ are the diagonal elements of the triangular matrix $L + D$, $i \in \{1,2,3\}$):
	\begin{equation}\label{40}
		|\frac{b}{a}|	< 1
	\end{equation}
	and the boundary of the convergence range~(\ref{39}), assuming that at least one of the roots of the equation~(\ref{38}) has a unit absolute value on it (let $r_1 = 1$), is given by the conditions:
	\begin{equation}\label{41}
		\begin{cases}
			d  = -ae^{i\varphi_{1}} -  b e^{-i\varphi_{1}}\\
			ar_{2}e^{i\varphi_{1}}e^{i\varphi_{2}} = b\\
			0\leq r_2 \leq 1
		\end{cases}
	\end{equation}
	\indent For $d_1 = \frac{d}{a}$ and $b_1 = \frac{b}{a}$ on the boundary~(\ref{41}), we can also find the relationships between absolute values and arguments:
	\[
	d_1 = -e^{i\varphi_{1}} - b_1 e^{-i\varphi_{1}}\\
	\]
	\[
	d_1 = -\cos\varphi_{1} - r_{b_1} \cos(\varphi_{b_1} - \varphi_{1}) - i(\sin\varphi_{1} + r_{b_1} \sin(\varphi_{b_1} - \varphi_{1}))
	\]
	\[
	\varphi_{b_1} = \varphi_{1} + \varphi_{2}
	\]
	\begin{equation}\label{42}
		|d_1| = \sqrt{1+r_{b_1}^2 + 2r_{b_1}\cos(\varphi_{1} - \varphi_{2})}
	\end{equation}
	\begin{equation}\label{43}
		\arg(d_1) = \arctan(\frac{\sin\varphi_{1} + r_{b_1}\sin\varphi_{2}}{\cos\varphi_{1} + r_{b_1}\cos\varphi_{2}})
	\end{equation}
	\indent In particular, when the matrix $A$ of the system~(\ref{1}) contains \textbf{real} coefficients, $a, b, d \in \mathbb{R}$, solving directly the quadratic equation~(\ref{38}) and applying the convergence criterion of the Gauss-Seidel method, taking into account the condition~(\ref{40}), we obtain the convergence range of the Gauss-Seidel method in the case of real roots of the equation~(\ref{38}):
	\begin{equation}\label{44}
		\begin{cases}
			|d|<|a+b| \\
			|\frac{b}{a}|<1
		\end{cases}
	\end{equation}
	and in the case of complex-conjugate roots of the equation~(\ref{38}).:
	\begin{equation}\label{45}
		0<\frac{b}{a}<1
	\end{equation}
	\indent Note that in the latter case, the condition
	\[
	|d|<|a+b|
	\]
	follows directly from the condition~(\ref{45}) and the negativity of the discriminant of the equation~(\ref{38}): $d^2 < 4ab$. \textbf{Therefore, the system~(\ref{44}) is a single range of convergence of the Gauss-Seidel method in the case of real matrix elements of the system~(\ref{1}).\\}
	\indent The first condition of the system~(\ref{44}) is interpreted as a segment $d$ on an infinite line.\\
	\indent Note that the conditions~(\ref{44}) are consistent with the boundary~(\ref{41}).\\
	\indent Unlike the Jacobi method, the convergence range of the Gauss-Seidel method in the case of real coefficients of the equation~(\ref{38}) is not constant, and the length of the above segment can vary depending on the parameters $a$ and $b$.\\
	\indent Let's compare the convergence ranges of both methods in the case of real coefficients of the system~(\ref{1}). For this purpose, we construct the convergence range bounded by the boundary~(\ref{34}) and the range~(\ref{44}) on the same coordinate plane $qOp$. The parameters $p$ and $q$ for the Jacobi method and $d, a, b$ for the Gauss-Seidel method are related by the relation:
	\[
	d = (p+q)a - b,
	\]
	substituting it into~(\ref{44}), we obtain:
	\begin{equation}\label{46}
		\begin{cases}
			|(p+q)a-b|<|a+b|\\
			|\frac{b}{a}|<1
		\end{cases}
	\end{equation}
	\indent Expanding the absolute values in the first inequality of the system~(\ref{46}), we find that one of the boundaries of the convergence range of the Gauss-Seidel method is always the line
	\begin{equation}\label{47}
		p=-q-1
	\end{equation}
	which is also one of the boundaries~(\ref{34}) of the convergence range of the Jacobi method, and the second one is also a straight line, which has the following form:
	\[
	p = -q + \frac{a+2b}{a}
	\]
	\indent It also shows that the size of the convergence range of the Gauss-Seidel method depends on the parameters $a, b$. Moreover, at some values of these parameters the convergence range of the Gauss-Seidel method can partially pass through the convergence range of the Jacobi method, and at other values it can completely contain it.\\
	\indent From the second inequality of the system~(\ref{46}) follows:
	\[
	-1<\frac{a+2b}{a}<3,
	\]
	so the convergence range of the Gauss-Seidel method on the plane $qOp$ is a part of this plane, which is always bounded from below by the line~(\ref{47}), and, depending on the particular case, bounded from above by a line parallel to it, the uppermost of which is the line
	\[
	p=-q+3.
	\]
	\newpage
	\indent Thus, together the convergence ranges of each method on the same plane $qOp$ are as follows (fig. \ref{p15}):
	\begin{figure}[h]
		\begin{center}
			\includegraphics[height=9cm, width=9.5cm]{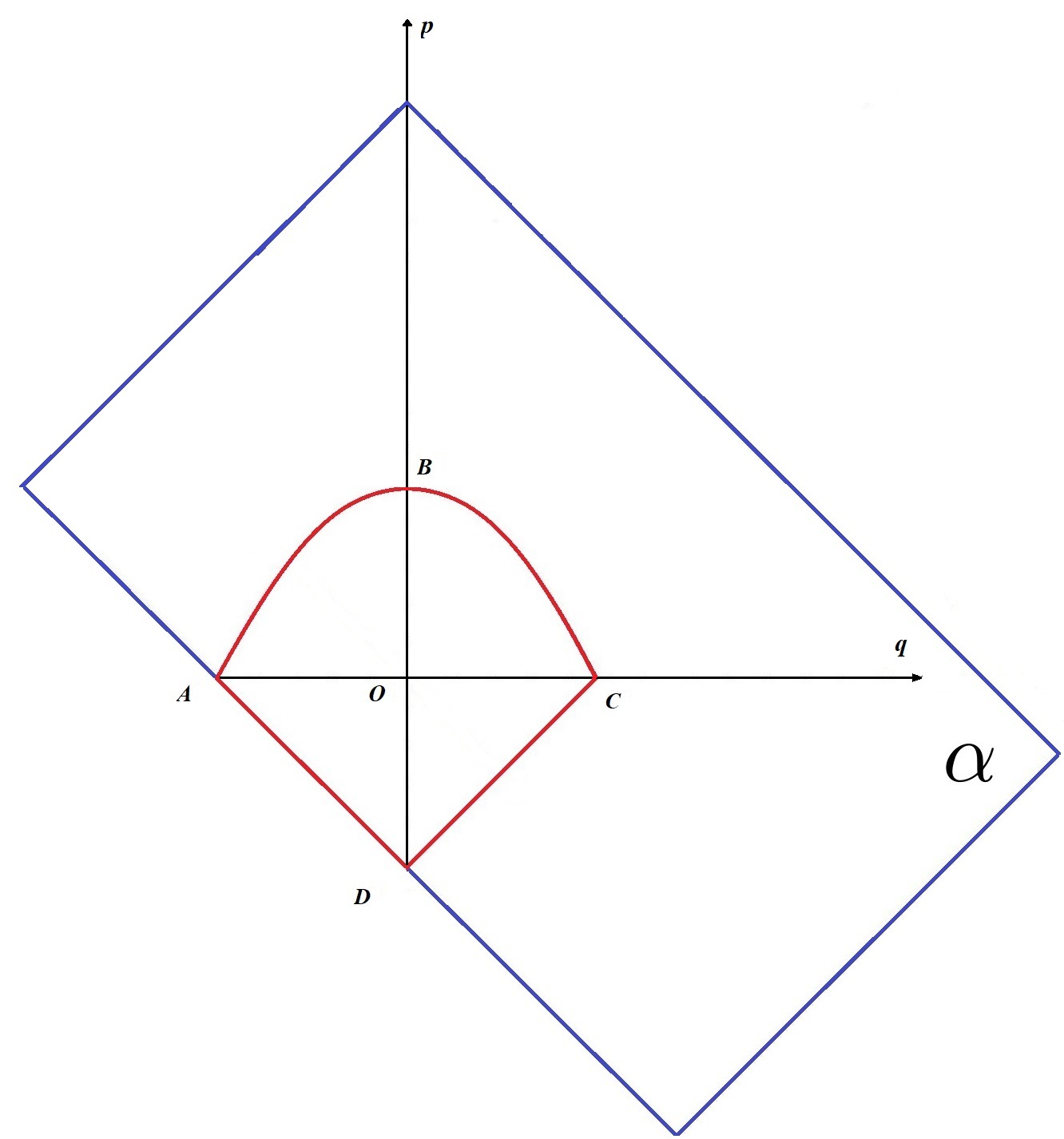}
			\caption{Convergence ranges of methods on the plane $qOp$. The band $\alpha$ is the maximum (with upper boundary $p=-q+3$) convergence range of the Gauss-Seidel method; the area $ABCD$ is the convergence range of the Jacobi method}
			\label{p15}
		\end{center}
	\end{figure}\\
	\indent According to figure \ref{p15}, the advantages of the Gauss-Seidel method over the Jacobi method when the system~(\ref{1}) has real matrix elements are obvious (in the case in Figure \ref{p15}, the convergence range of the Jacobi method is entirely contained in the convergence range of the Gauss-Seidel method), especially when the parameters $p$ and $q$ have large absolute values - then the Jacobi method does not converge. Nevertheless, the upper bound of the range for the Gauss-Seidel method varies depending on the parameters $a$ and $b$, so if the iterative process of the Jacobi method converges for the SLAE, it does not mean that the iterative process of the Gauss-Seidel method converges.\\
	\indent Let's give examples of constructing the convergence range of the Gauss-Seidel method in coordinates $qOp$ to demonstrate how it varies depending on the parameters $a$ and $b$, and in the same coordinates we construct the convergence range of the Jacobi method for clarity.\\
	\indent \textbf{Example 1.}\\
	\indent Let the parameters $a=2$, $b=1$, then the convergence range of the Gauss-Seidel method has the form:
	\[
	\begin{cases}
		|2(p+q)-1|<3\\
		|\frac{b}{a}|= \frac{1}{2} < 1
	\end{cases}
	\]
	thus 
	\[
	\frac{a+2b}{a}=2
	\]
	\indent Then, by analogy with fig. ~\ref{p15}, the convergence ranges for each method on the plane $qOp$ look as follows (fig. ~\ref{p16}):
	\begin{figure}[h]
		\center
		\includegraphics[height=6cm, width=11cm]{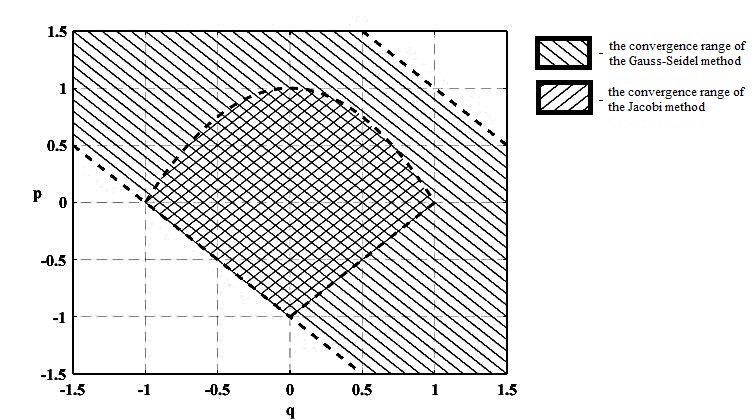}
		\caption{Convergence ranges of Jacobi and Gauss-Seidel methods at parameters $a = 2$ and $b = 1$}
		\label{p16}
	\end{figure}\\
	\indent The figure~\ref{p16} shows that the convergence range of the Jacobi method lies entirely within the convergence range of the Gauss-Seidel method, so in this particular case of parameters $a, b$ for any SLAE for which the Jacobi method converges, the Gauss-Seidel method also converges, but the converse is not true.\\
	\indent \textbf{Example 2.}\\
	\indent Here is an example of a SLAE in three unknowns, for which the Jacobi method converges, but the Gauss-Seidel method does not converge:
	\[
	A =
	\begin{pmatrix}
		-8& 6&-4\\
		-9& 8& 6\\
		4&-5& 3
	\end{pmatrix}
	\]
	\\
	\indent In this case the parameters are as follows: $a = -192$, $b = 144$, then the convergence range of the Gauss-Seidel method has the following form:
	\[
	\begin{cases}
		|-192(p+q)-144|<48\\
		|\frac{b}{a}|=|\frac{144}{192}|<1
	\end{cases}
	\]
	\[
	\frac{a+2b}{a}= -\frac{1}{2}
	\]
	\indent By analogy with fig. ~\ref{p15}, we obtain the following picture of convergence ranges for both methods on the plane $qOp$ (fig. ~\ref{p17}):
	\begin{figure}[h]
		\center
		\includegraphics[height=6cm, width=11cm]{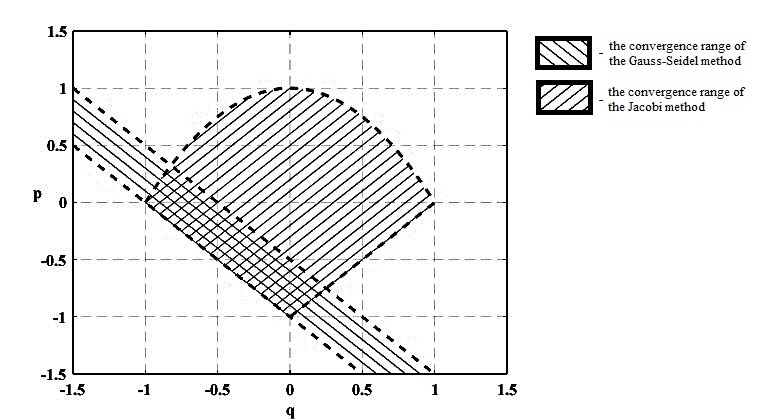}
		\caption{Convergence ranges of Jacobi and Gauss-Seidel methods at parameters $a = -192$ and $b = 144$}
		\label{p17}
	\end{figure}\\
	\newpage
	\indent This range does not satisfy the parameters $p$, $q$, which in this particular example for matrix $A$ are equal to:
	\[
	p=-\frac{50}{192} 
	\]
	\[
	q= \frac{36}{192} 
	\]
	\indent From the figure~\ref{p17} it is obvious that this point ($p,q$) does not belong to the convergence range of the Gauss-Seidel method on the plane $qOp$, but it belongs to the convergence range of the Jacobi method.\\
	\indent In addition, in this particular case we see that the Gauss-Seidel method does not converge in most of the convergence range of the Jacobi method, but it can converge at large values of $p$ and $q$, while the Jacobi method does not converge at large values of $p$ and $q$.
	
	\section{The general case of systems of linear algebraic equations with complex matrices}
	\label{GenCase} 
	\indent The convergence check of each method is an investigation to find all roots of a polynomial of degree $n$ inside the unit circle, which can be transformed to a stability study problem \cite{Zadorojniy}.\\
	\indent In general, a polynomial of degree $n$ with complex coefficients is obtained from the determinant equations~(\ref{2}) or~(\ref{3}):
	\begin{equation}\label{48}
		f(\lambda)=a_0\lambda^n + a_1\lambda^{n-1} + ... + a_n = 0, a_0 \neq 0
	\end{equation}
	\indent For convergence of the method to which the given polynomial corresponds, it is necessary and sufficient that all its roots lie inside the unit circle, for which, in turn, it is necessary and sufficient that the polynomial:
	\[
	f(z)=a_0(z+1)^n + a_1(z+1)^{n-1}(z-1) + a_2(z+1)^{n-2}(z-1)^2 +...+ a_n(z-1)^n = 0
	\]
	obtained from~(\ref{48}) be stable \cite{Zadorojniy}.\\
	\indent In general, to check its stability, we can use the complex analog of Hurwitz's stability criterion \cite{Postnikov}: let there be an arbitrary polynomial of degree $n$ with complex coefficients, the stability of which should be investigated:
	\[
	f(z)=d_0z^n + d_1z^{n-1} + ... + d_n,
	\]
	 it's equivalent, under the assumption that $d_0 \neq 0$, to the polynomial whose first coefficient is equal to one:
	\begin{equation}\label{49}
		\widetilde{f}(z)=z^n + \frac{d_1}{d_0}z^{n-1} + ... + \frac{d_n}{d_0}
	\end{equation}
	\indent Replacing in~(\ref{49}) the variable $z$ by a purely imaginary number $i\omega, \omega \in \mathbb{R}$, we have the polynomial:
	\[
	\widetilde{f}(i\omega)=(i\omega)^n + \frac{d_1}{d_0}(i\omega)^{n-1} + ... + \frac{d_n}{d_0}
	\]
	which, by raising the multiplier $i\omega$ of each summand to the appropriate degree and separating the purely imaginary elements from the purely real ones, is represented as the sum of two polynomials with real coefficients.
	\[
	\widetilde{f}(i\omega)=g(\omega) + ih(\omega),
	\]
	 for which, according to \cite{Postnikov}, if the degree of the polynomial~(\ref{49}) is $n = 2m$:
	\[
	\widetilde{g} = (-1)^m g
	\]
	\[
	\widetilde{h} = (-1)^{m-1} h,
	\]
	if $n = 2m + 1$:
	\[
	\widetilde{g} = (-1)^m h
	\]
	\[
	\widetilde{h} = (-1)^m g.
	\]
	\indent Let 
	\[
	B = b_0x^n + b_1x^{n-1} +...+b_n
	\]
	be an arbitrary polynomial of degree $n$ with real coefficients with positive prime factor $b_0$, and let
	\[
	C = c_0x^{n-1}  + c_1x^{n-2} + ... + c_{n-1}
	\]
	be an arbitrary polynomial of degree at most $n-1$ with real coefficients.\\
	\indent \textbf{Definition.} Square matrix 
	\[
	\begin{pmatrix}
		b_{0}& b_{1}& b_{2}&b_{3}& ...& 0\\
		0& c_{0}& c_{1}&c_{2}& ...& 0 \\
		0& b_{0}& b_{1}&b_{2}& ...& 0 \\
		0& 0& c_{0}&c_{1}& ...& 0 \\
		...&...&...&...&...&...\\
		0&...&b_{0}&b_{1}& ...& b_{n} \\
		0&...& 0   & c_{0}&...&c_{n-1}
	\end{pmatrix}
	\]
	of order $2n$ is called the Hurwitz matrix of polynomials $B$ and $C$, and its principal minors of \textit{even order} are called \textbf{the Hurwitz determinants} of polynomials $B$ and $C$.\\
	\indent \textbf{The complex analog of the Hurwitz stability criterion:} polynomial of degree $n$
	\[
	\widetilde{f}(z)=z^n + \frac{d_1}{d_0}z^{n-1} + ... + \frac{d_n}{d_0} = 0
	\]
	with complex coefficients and a unit (real and positive, but not necessarily unit) coefficient at the highest degree is stable if and only if all Hurwitz determinants of polynomials $\widetilde{g}$ and $\widetilde{h}$ are positive. \\
	\indent In particular, when all coefficients of the resulting polynomial are real, we can use the classical Rouse-Hurwitz stability criterion for polynomials with real coefficients, or other similar \cite{Gantmaher} criteria to check stability.\\
	\indent Thus, in general, to check the convergence of the Jacobi and Gauss-Seidel iterative methods, in order to avoid a direct search for the roots of a polynomial with complex or real coefficients, it is necessary to reduce it to a new one, which is checked for stability, which can be done using a computer by the above method. This method of checking convergence is especially relevant when the initial SLAEs have a large dimension, because of which we obtain equations of large powers, the solution of which is often very cumbersome. \\
	\indent Let's show that the range bounded by the boundary~(\ref{34}) is also obtained by applying the described method of checking convergence through the complex analog of the Hurwitz criterion for a polynomial with complex coefficients:
	\[
	f(\lambda) = \lambda^3 + p\lambda + q
	\]
	\[
	k(z) = (z+1)^3 + p(z+1)(z-1)^2 + q(z-1)^3
	\]
	\[
	k(z) = (1+p+q)z^3 + (3-p-3q)z^2 + (3-p+3q)z + (1+p-q)
	\]
	\indent Assuming that $1+p+q \neq 0$, divide the last polynomial by this sum
	\[
	\widetilde{k}(z) = z^3 + \frac{3-p-3q}{1+p+q}z^2 + \frac{3-p+3q}{1+p+q}z + \frac{1+p-q}{1+p+q}
	\]
	\[
	\widetilde{k}(iw) = -iw^3 - \frac{3-p-3q}{1+p+q} w^2 + \frac{3-p+3q}{1+p+q}iw + \frac{1+p-q}{1+p+q}
	\]
	\indent Let's separate the real and imaginary parts
	\[
	\widetilde{k}(iw) = (-\operatorname{Re}\frac{3-p-3q}{1+p+q} w^2 - \operatorname{Im}\frac{3-p+3q}{1+p+q}w + \operatorname{Re}\frac{1+p-q}{1+p+q}) +
	\]
	\[
	+ i(-w^3 - \operatorname{Im}\frac{3-p-3q}{1+p+q}w^2 + \operatorname{Re}\frac{3-p+3q}{1+p+q}w + \operatorname{Im}\frac{1+p-q}{1+p+q})
	\]
	\[
	g(w) = -\operatorname{Re}\frac{3-p-3q}{1+p+q} w^2 - \operatorname{Im}\frac{3-p+3q}{1+p+q}w + \operatorname{Re}\frac{1+p-q}{1+p+q}
	\]
	\[
	h(w) = -w^3 - \operatorname{Im}\frac{3-p-3q}{1+p+q}w^2 + \operatorname{Re}\frac{3-p+3q}{1+p+q}w + \operatorname{Im}\frac{1+p-q}{1+p+q}
	\]
	\indent The degree of the polynomial $\widetilde{k}(z)$ is odd, so
	\[
	\tilde{g} = -h(w) = w^3 + \operatorname{Im}\frac{3-p-3q}{1+p+q}w^2 - \operatorname{Re}\frac{3-p+3q}{1+p+q}w - \operatorname{Im}\frac{1+p-q}{1+p+q}
	\]
	\[
	\tilde{h} = -g(w) = \operatorname{Re}\frac{3-p-3q}{1+p+q} w^2 + \operatorname{Im}\frac{3-p+3q}{1+p+q}w - \operatorname{Re}\frac{1+p-q}{1+p+q}
	\]
	\indent The Hurwitz matrix for polynomials $\tilde{g}, \tilde{h}$ has the form:
	\[
	\begin{pmatrix}
		1& \operatorname{Im}\frac{3-p-3q}{1+p+q}& -\operatorname{Re}\frac{3-p+3q}{1+p+q}& -\operatorname{Im}\frac{1+p-q}{1+p+q}&0&0\\
		0& \operatorname{Re}\frac{3-p-3q}{1+p+q}& \operatorname{Im}\frac{3-p+3q}{1+p+q}& -\operatorname{Re}\frac{1+p-q}{1+p+q}&0&0\\
		0&1&\operatorname{Im}\frac{3-p-3q}{1+p+q}&-\operatorname{Re}\frac{3-p+3q}{1+p+q}&-\operatorname{Im}\frac{1+p-q}{1+p+q}&0\\
		0&0&\operatorname{Re}\frac{3-p-3q}{1+p+q}&\operatorname{Im}\frac{3-p+3q}{1+p+q}&-\operatorname{Re}\frac{1+p-q}{1+p+q}&0\\
		0&0&1& \operatorname{Im}\frac{3-p-3q}{1+p+q}& -\operatorname{Re}\frac{3-p+3q}{1+p+q}& -\operatorname{Im}\frac{1+p-q}{1+p+q}\\
		0&0&0&\operatorname{Re}\frac{3-p-3q}{1+p+q}& \operatorname{Im}\frac{3-p+3q}{1+p+q}& -\operatorname{Re}\frac{1+p-q}{1+p+q}
	\end{pmatrix}
	\]
	\indent When finding the convergence range in the real case, equating all imaginary elements in the obtained Hurwitz matrix to zero, we obtain the corresponding Hurwitz matrix, the principal minors of even order of which give the conditions we obtained above from the boundary~(\ref{34}).\\
	\indent In general, to check the convergence of the Jacobi method for a particular SLAE in three unknowns (in our case), we can program the described algorithm. For example, in the Python language:
	\begin{figure}[h]
		\centering
		\includegraphics[width=\textwidth]{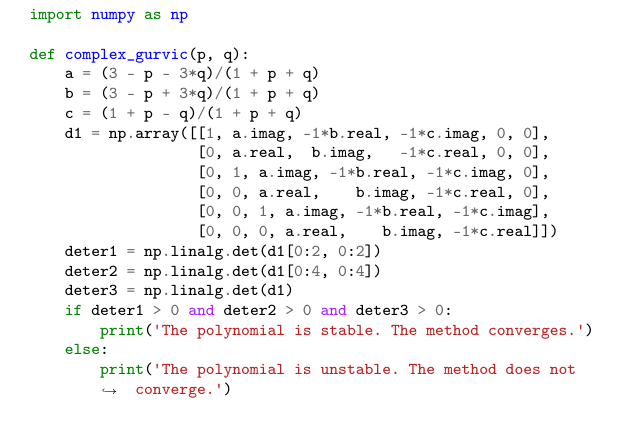}
	\end{figure}
	
	\indent In general, the following conclusion can be made about the comparison of convergence of the two methods: in the equation~(\ref{3}) for the Gauss-Seidel method, it is always possible to take $\lambda$ from the last line beyond the sign of the determinant, thus lowering the degree of the polynomial whose stability is to be investigated by one, which is not always possible for the Jacobi method according to the equation~(\ref{2}). Thus, for SLAEs in $n>2$ unknowns with complex matrices, in general case, the polynomial, whose stability should be investigated, obtained for the Jacobi method, has degree by one more in contrast to the analogous polynomial for the Gauss-Seidel method.
	
	\section{Statistical comparison of convergence of Jacobi and Gauss-Seidel methods}
	\label{Stat} 
	\indent 100000 random matrices of SLAEs~(\ref{1}) with real matrix elements that are uniformly distributed random variables on the interval $[-100; 100]$, with the number of unknowns from two to five, for each of them the well-known convergence criteria of each method were checked, then for each number of unknowns the number of cases in which both methods converge, only the Gauss-Seidel method converges, only the Jacobi method converges was determined. The obtained data are summarized in the table~\ref{t1}.
	
	\begin{table}[h]
		\caption{Convergence results of Jacobi and Gauss-Seidel methods}
		\begin{tabular}{ | c | c | c |c|}
			\hline
			\thead{Number of \\ unknowns} & \thead{Both methods \\converge} & \thead{The Gauss-Seidel method\\ converges, but the Jacobi method \\does not converge}& \thead{The Jacobi method converges,\\ but the Gauss-Seidel method \\does not converge} \\ 
			\hline
			2&49916&0&0\\ 
			\hline
			3&11818&7521&1095\\ 
			\hline
			4&1436&3411&528 \\ 
			\hline
			5 & 111 & 726 & 76 \\ 
			\hline
		\end{tabular}
		\label{t1}
	\end{table}
	\indent The data obtained in the table for the number of unknowns $n>2$ confirm the conclusions that, in general, the Gauss-Seidel method converges much more often than the Jacobi method, but the convergence of one of the methods cannot guarantee the convergence of the other. At the same time, we also see that as the number of unknowns in the SLAEs increases, both methods converge much less frequently, which is consistent with the above complex analog of the Hurwitz criterion.\\
	\indent Note also that in the case of SLAEs in two unknowns, the data from the table~\ref{t1} confirm the conclusions that in this case both methods converge in the same way - if one converges, the other converges as well.
	
	\section{Conclusion}
	\label{Conc} 
	\indent The found boundary conditions in the complex case, as well as convergence ranges in the real case allowed us to see the picture of convergence conditions of Jacobi and Gauss-Seidel iterative methods and on this basis to make a comparative analysis of the effectiveness of each method: if in the case of square matrices of SLAEs in two unknowns both methods converge equally effectively, in the case of matrices of SLAEs in three and more unknowns methods have a noticeable difference in the convergence conditions - with increasing number of unknowns in SLAEs, the Gauss-Seidel method is noticeably more effective.\\
	\indent For example, in the case of an SLAEs' matrices in three unknowns, when the convergence ranges for both methods are plotted for the real case on the same coordinate plane, it can be seen that in the general case the Gauss-Seidel method has better convergence than the Jacobi method, since its convergence range is bounded by straight lines, but infinite in contrast to the convergence range of the Jacobi method, one of whose boundaries even enters the boundary of the convergence range of the Gauss-Seidel method. However, as it has been shown, the convergence range of the Gauss-Seidel method depends on the parameters that do not always give a full convergence range of the Jacobi method into the convergence range of the Gauss-Seidel method, because of which there may be situations when iterations converge to the exact solution by the Jacobi method but do not converge by the Gauss-Seidel method. Statistical comparison of convergence of both methods also confirms these conclusions.\\
	\indent When the number of unknowns over the field of complex numbers is large, the convergence of each method can be checked using the complex analog of the Hurwitz stability criterion, or using the classical Rouse-Gurwitz criterion in the real case.\\

\end{document}